\documentclass{amsart}

\usepackage{amsmath,amssymb,amsfonts}

\newtheorem{theorem}{Theorem}
\newcommand{\bt}{\begin{theorem}}
\newcommand{\et}{\end{theorem}}
\newtheorem{lemma}{Lemma}
\newcommand{\bl}{\begin{lemma}}
\newcommand{\el}{\end{lemma}}
\newcommand{\beq}{\begin{equation}}
\newcommand{\eeq}{\end{equation}}
\newcommand{\mca}{\ensuremath{ \mathcal A}}
\newcommand{\mcf}{\ensuremath{ \mathcal F}}
\newcommand{\N}{\ensuremath{ \mathbf N }}
\newcommand{\Z}{\ensuremath{\mathbf Z}}
\newcommand{\benum}{\begin{enumerate}}
\newcommand{\eenum}{\end{enumerate}}

\title{Shatrovski\u{i}'s construction of thin bases}
\author{Melvyn B. Nathanson}
\address{Department of Mathematics\\
Lehman College (CUNY), Bronx, NY 10468\\
and CUNY Graduate Center, New York, NY 10016}
\email{melvyn.nathanson@lehman.cuny.edu}

\thanks{Supported in part by a grant from the PSC-CUNY Research Award Program.}

\date{\today}

\begin{document}
\maketitle

\begin{abstract}
The set \mca\ of nonnegative integers is called a basis of order $h$ if every nonnegative integer can be represented as the sum 
of exactly $h$ not necessarily distinct elements of \mca.  
An additive basis \mca\ of order $h$ is called thin if there exists $c > 0$ such that the number of elements of $\mca$ not exceeding $x$ is less than $c x^{1/h}$ for all $x \geq 1$.  
This paper describes a construction of Shatrovski\u{i} of thin bases of order $h$.
\end{abstract}

\section{Additive bases of finite order}
Let $\N_0$ denote the set of nonnegative integers, let $h \geq 2$, and let $(\mca_1,\ldots, \mca_h)$ be an $h$-tuple of subsets of $\N_0$.  We define the \emph{sumset}
\[
\mca_1 + \cdots + \mca_h = \{a_1 + \cdots+ a_h : a_i \in \mca_i \text{ for } i = 1,\ldots, h\}.
\]
The \emph{$h$-fold sumset} of a set $\mca$ of nonnegative integers is 
\[
h\mca = \underbrace{\mca+ \cdots + \mca}_{h \text{ terms}} .
\]
The $h$-tuple $(\mca_1,\ldots, \mca_h)$ is called an \emph{additive system} of order $h$ if $\mca_1 + \cdots + \mca_h = \N_0.$
The set \mca\ is an \emph{additive basis}, or, simply, a \emph{basis},  of order $h$ if $h\mca = \N_0.$  If $(\mca_1,\ldots, \mca_h)$ is  an additive system of order $h$, then $\mca = \bigcup_{i=1}^h \mca_i$ is a basis of order $h$.  

For example, let $a^{\prime}_1,\ldots, a^{\prime}_h$ be relatively prime positive integers, and let
\[
\mca_i = a^{\prime}_i\ast \N_0 = \{a^{\prime}_iv_i : v_i \in \N_0\}
\]
for $i=1,\ldots, h$.  
There exists an integer $C$ such that every integer $n \geq C$ can be represented in the form $n = a^{\prime}_1v_1 + \cdots + a^{\prime}_hv_h$ with $v_1,\ldots, v_h \in \N_0$ (Nathanson~\cite[Section 1.6]{nath00aa}, and so the set
\[
\mca = [0,C-1] \cup \bigcup_{i=1}^h \mca_i
\]
is an additive basis of order $h$.

Let $g \geq 2$.  We can use the $g$-adic representation of the nonnegative integers to construct another class of additive bases of order $h$.  Let  $K_1,\ldots, K_h$ be subsets of $\N_0$ such that 
\[
\N_0 = K_1\cup \cdots \cup K_h.
\]
For $i = 1,\ldots, h$, let $\mcf(K_i)$ denote the set of all finite subsets of $K_i$, and let 
\[
\mca_i = \left\{ 
\sum_{k\in F} \varepsilon_k g^k : \varepsilon_k \in \{ 0,1,2,\ldots, g-1\} \text{ and } F \in \mcf(K_i)  
\right\}.
\]
Then $(\mca_1,\ldots, \mca_h)$ is  an additive system of order $h$ and $\mca = \bigcup_{i=1}^h \mca_i$ is a basis of order $h$.  

The \emph{counting function} of a set \mca\ of integers is
\[
A(x) = \sum_{\substack{a\in\mca\\1 \leq a \leq x}} 1.
\]
If \mca\ is a basis of order $h$, then for all $x \in \N_0$ we have
\[
x+1 \leq {A(x)+h \choose h} \leq \frac{ (A(x)+h)^h}{h!}
\]
and so 
\[
\liminf_{x\rightarrow\infty} \frac{A(x)}{x^{1/h}} \geq (h!)^{1/h}.
\]
The basis \mca\  of order $h$ is called a \emph{thin basis} if  
\[
\limsup_{x\rightarrow\infty} \frac{A(x)}{x^{1/h}} < \infty.
\]
One of the first problems in additive number theory for general sets of integers concerned the existence and construction of thin bases of finite order (Rohrbach~\cite{rohr37a}).    
In 1937, Raikov~\cite{raik37} and St\" ohr~\cite{stoh37}  independently solved this problem. Using the 2-adic representation and the sets
$K_i = \{ k\in \N_0 : k \equiv i-1 \pmod{h} \}$ for $i=1,\ldots, h$, they proved that $\mca = \bigcup_{i=1}^h \mca_i$ is a thin basis of order $h$. 
Many thin bases have been constructed using variations of the $g$-adic construction, for example, Cassels~\cite{cass57,grek-hadd-helo-pihk06,nath09xf} and Jia-Nathanson~\cite{nath89b}.  
Recently,  Hofmeister~\cite{hofm01}, Blomer~\cite{blom03}, and Schmitt~\cite{schm06} constructed other examples of thin bases.  

In 1940 Shatrovski\u{i}~\cite{shat40} described a beautiful class of thin bases of order $h$ that use the linear diophantine equation $n = a^{\prime}_1x_1+\cdots a^{\prime}_hx_h$ and do not depend on $g$-adic representations.   This is an almost  forgotten paper.\footnote{Shatrovski\u{i} published his paper in Russian, with a summary in French.  The transliteration of the author's name in the French summary is Chatrovsky, and this is the form of the name that appears in  \emph{Mathematical  Reviews}.  I use the AMS and Library of Congress system for transliteration of Cyrillic into English.}
 It has been mentioned in a few surveys of solved and unsolved problems in additive number theory (e.g. St\" ohr~\cite{stoh55}, Erd\H os and Nathanson~\cite{nath87e}, and Nathanson~\cite{nath89e}), but  MathSciNet lists no publication that cites Shatrovski\u{i}'s paper.
This note provides an exposition and generalization of Shatrovski\u{i}'s construction of thin bases.

\section{The construction}
Fix an integer $h \geq 2$,  and define the positive integer 
\beq    \label{thin:shat-k0}
k_0 = \left[ \frac{1}{2^{1/h}-1} \right] + 1.
\eeq
Then  $1 + (1/k) < 2^{1/h}$ for all $k \geq k_0.$  
Let $r_1,\ldots, r_h$ be pairwise relatively prime positive integers, and arrange them in increasing order so that $r_i < r_{i+1}$ for $i = 1,\ldots, h-1$.   Let $P$ be a positive integer such that (i) $P \geq r_h - r_1$, and (ii) if $1 \leq i < j \leq h$ and if the prime $p$ divides $r_j-r_i$, then $p$ divides $P$.  
For every positive integer $k$, we define 
\beq    \label{thin:shat-sik}
s_{i,k} = kP+r_i
\eeq
for $i = 1,\ldots,h$.    The inequality $r_1 < r_2 < \cdots < r_h$ implies  that, for all $k \geq 1$, 
\[
s_{1,k} < s_{2,k} < \cdots < s_{h,k}.
\]
We define the integers 
\beq    \label{thin:shat-Sk}
S_k = \prod_{i=1}^h s_{i,k}
\eeq
and
\beq    \label{thin:shat-aik}
a_{i,k} = \frac{S_k}{s_{i,k}} = \prod_{\substack{ j=1 \\ j\neq i}}^h s_{j,k}
\eeq
for $i=1,\ldots, h$.   For all $k \geq 1$ we have $S_k > 2^h$ and
\[
a_{1,k} > a_{2,k} > \cdots > a_{h,k}.
\]

\bl     \label{thin:lemma:shat1}
Define $k_0$ by~\eqref{thin:shat-k0}.  Let $r_1,\ldots, r_h$ be a strictly increasing sequence of pairwise relatively prime positive integers, and let $s_{1,k},\ldots, s_{h,k}$, $S_k$, and $a_{1,k},\ldots, a_{h,k}$ be the positive integers defined by~\eqref{thin:shat-sik},~\eqref{thin:shat-Sk}, and~\eqref{thin:shat-aik}, respectively.  
For all $k \geq 1$,
\benum
\item[(i)]
the integers $s_{1,k},  \ldots, s_{h,k}$ are pairwise relatively prime,  
\item[(ii)]
the integers $a_{1,k},  \ldots, a_{h,k}$ are relatively prime, 
\item[(iii)]
 $S_k < S_{k+1}$ and, if $k \geq k_0$, then $S_{k+1} < 2S_k.$
\eenum
\el

\begin{proof}
Let $1 \leq i < j \leq h$.  If $(s_{i,k},s_{j,k}) = d >1$, then $d$ divides both $s_{i,k}$ and $s_{j,k}$, and so $d$ divides 
$s_{j,k} - s_{i,k} = r_j - r_i $.  If $p$ is a prime divisor of $d$, then $p$ divides $r_j - r_i $ and so $p$ divides $P$.  Therefore, $p$ divides both $r_i$ and $r_j$, which is absurd, since $(r_i,r_j) = 1$.  We conclude that $s_{1,k}, s_{2,k}, \ldots, s_{h,k}$ are pairwise relatively prime.  This proves~(i).

If $\gcd( a_{1,k}, a_{2,k}, \ldots, a_{h,k} ) > 1$, then there is a prime $p$ that divides $a_{i,k}$ for all $i = 1,\ldots, h$.  Since $p$ divides $a_{1,k}$, it follows from~\eqref{thin:shat-aik} that $p$ divides $s_{j,k}$ for some $j \in \{2,3,\ldots, h\}$.  Similarly, since $p$ divides $a_{j,k}$, it follows from~\eqref{thin:shat-aik} that $p$ divides $s_{\ell,k}$ for some $\ell \neq j$.  This is impossible, since $(s_{j,k},s_{\ell,k}) =1$, and so $( a_{1,k}, a_{2,k}, \ldots, a_{h,k} ) = 1$.   This proves~(ii).

To prove~(iii), we use the inequality 
\[
s_{i,k} = kP+r_i < (k+1)P+r_i = s_{i,k+1}
\]
for $i = 1,\ldots, h$ to obtain $S_k < S_{k+1}$.  
Since $(x+z)/(y+z) < x/y$ if $x > y > 0$ and $z > 0$, it follows that for $k \geq k_0$ we have 
\[
\frac{s_{i,k+1}}{s_{i,k}} = \frac{(k+1)P + r_i}{k P+r_i} 
< \frac{(k+1)P }{kP } =1 +\frac{1}{k}  \leq 1 +\frac{1}{k_0} < 2^{1/h}
\]
and so
\[
\frac{S_{k+1}}{S_k} = \prod_{i=1}^h \frac{s_{i,k+1}}{s_{i,k}} < \left( 1 +\frac{1}{k_0} \right)^h < 2.
\]
This completes the proof.  
\end{proof}

\bl         \label{thin:lemma:shat2}
Let  $\{S_k\}_{k=0}^{\infty}$ be the strictly increasing sequence of positive integers defined by~\eqref{thin:shat-Sk}.   
Choose a positive integer $k_1 \geq k_0$. 
For every positive integer $t$, there is a unique integer $\ell_t$ such that 
\beq  \label{thin:ShatSeq}
S_{k_1+\ell_t} \leq S_{k_1}^{t} < S_{k_1+\ell_{t -1} +1}.
\eeq
Then $\ell_1 = 0$ and  the sequence $\{\ell_t\}_{t=1}^{\infty}$ is strictly increasing. 
\el

\begin{proof}
We have $S_{k_1} \leq S_{k_1}^1 < S_{k_1+1}$ and so $\ell_1 = 0$.   By Lemma~\ref{thin:lemma:shat1} (iii),  we have $S_{k+1} < 2S_k$ for all $k \geq k_0$.  Since $S_k > 2^h$ for all $k\geq 1$, it follows that 
\[
S_{k_1}^t < S_{k_1+\ell_t + 1} < 2S_{k_1+\ell_t} < 2^h S_{k_1}^t 
< S_{k_1}^{t+1} < S_{k_1 + \ell_{t+1} + 1}
\]
and so $\ell_t  <  \ell_{t+1}$.  This completes the proof.  
\end{proof}

\bl       \label{thin:lemma:shat3}
Let $r_1,\ldots, r_h$ be a strictly increasing sequence of pairwise relatively prime positive integers, let $k \geq 1$, and let $s_{1,k},\ldots, s_{h,k}$, $S_k$, and $a_{1,k},\ldots, a_{h,k}$ be the integers defined by~\eqref{thin:shat-sik},~\eqref{thin:shat-Sk}, and~\eqref{thin:shat-aik}, respectively.   If 
\[
(h-1)S_k < n \leq L
\]
then there exist nonnegative integers $v_1, v_2, \ldots, v_h$ such that 
\[
a_{1,k}v_1 + a_{2,k}v_2 + \cdots + a_{h,k}v_h = n
\]
and
\[
v_i \leq  \frac{L}{a_{i,k}}
\]
for $i=1,\ldots, h.$
\el

\begin{proof}
Since $\gcd( a_{1,k}, a_{2,k}, \ldots, a_{h,k} ) = 1$, there exist integers $u_1,u_2,\ldots,u_h$, not necessarily nonnegative, such that 
\[
a_{1,k}u_1 + a_{2,k}u_2 + \cdots + a_{h,k}u_h = n.
\]
For each $i = 1,2,\ldots, h-1$ there is a unique integer $v_i \in \{0,1,2,\ldots, s_{i,k}-1\}$  such that 
\[
v_i \equiv u_i \pmod{s_{i,k}}.
\]
Since $S_k = a_{i,k}s_{i,k}$ for $i = 1,\ldots, h$, we have 
\[
 a_{1,k}v_1 + a_{2,k}v_2 + \cdots + a_{h-1,k}v_{h-1} \leq 
 \sum_{i=1}^{h-1}a_{i,k}(s_{i,k}-1) <  (h-1)S_k
 \]
and
\begin{align*}
n & =  a_{1,k}u_1 + a_{2,k}u_2 + \cdots + a_{h-1,k}u_{h-1} + a_{h,k}u_h \\
& =  a_{1,k}v_1 + a_{2,k}v_2 + \cdots + a_{h-1,k}v_{h-1} + a_{h,k}u_h + \sum_{i=1}^{h-1} a_{i,k} (u_i-v_i) \\
& = a_{1,k}v_1 + a_{2,k}v_2 + \cdots + a_{h-1,k}v_{h-1} + 
a_{h,k}u_h + S_k \sum_{i=1}^{h-1} \left( \frac{u_i-v_i}{s_{i,k}} \right)  \\
& = a_{1,k}v_1 + a_{2,k}v_2 + \cdots + a_{h-1,k}v_{h-1} + 
a_{h,k}u_h + a_{h,k}s_{h,k} \sum_{i=1}^{h-1} \left( \frac{u_i-v_i}{s_{i,k}} \right)  \\
& = a_{1,k}v_1 + a_{2,k}v_2 + \cdots + a_{h-1,k}v_{h-1} + a_{h,k}v_h
\end{align*}
where 
\[
v_h = u_h + s_{h,k} \sum_{i=1}^{h-1} \left( \frac{u_i-v_i}{s_{i,k}} \right) \in \Z.
\]
If $n \geq (h-1)S_k$, then 
\[
a_{h,k}v_h = n - \sum_{i=1}^{h-1} a_iv_i >  n - (h-1)S_k \geq 0
\]
and so $v_h \geq 1$.  
For $i=1,\ldots, h$, the inequality  
\[
0 \leq a_{i,k}v_i \leq n \leq L
\]
implies that  
\[
v_i \leq \frac{L}{a_{i,k}}.
\]
This completes the proof.
\end{proof}

\bt       \label{thin:theorem:basis}
Let $h \geq 2$ and $k_1 \geq k_0$.  
Let $r_1,\ldots, r_h$ be a strictly increasing sequence of pairwise relatively prime positive integers, and let $s_{1,k},\ldots, s_{h,k}$, $S_k$, and $a_{1,k},\ldots, a_{h,k}$ be the integers defined by~\eqref{thin:shat-sik},~\eqref{thin:shat-Sk}, and~\eqref{thin:shat-aik}, respectively.  Let $\{\ell_t\}_{t=0}^{\infty}$ be the increasing sequence of nonnegative integers defined by~\eqref{thin:ShatSeq}.  
For $t \geq 1$, let
\[
V(t) = \left\{a_{i,k_1+\ell_t} v_i  : i = 1,\ldots, h \text{ and }
1 \leq v \leq  \frac{(h-1)S_{k_1+ \ell_{t+1}}}{ a_{i,k+\ell_t} } \right\}.
\]
The set
\[
\mca = [0,(h-1)S_{k_1}] \cup \bigcup_{t = 1}^{\infty}   V(t)
\]
is a thin basis of order $h$.
\et

\begin{proof}
Since $0 \in \mca$, we have $[0,(h-1)S_{k_1}] \subseteq \mca \subseteq h\mca.$  Recall that $\ell_1 = 0$.  
By Lemmas~\ref{thin:lemma:shat1} and~\ref{thin:lemma:shat2}, the sequences 
$\{S_k\}_{k=k_1}^{\infty}$ and $\{\ell_t\}_{t=1}^{\infty}$ are strictly increasing.  Therefore, for every $n > (h-1)S_{k_1}$ there is a unique positive integer $t$ such that 
\[
(h-1)S_{{k_1}+ \ell_t} < n \leq (h-1)S_{{k_1}+ \ell_{t+1}}.
\]
By Lemma~\ref{thin:lemma:shat3}, the linear diophantine equation
\[
a_{1,{k_1}+ \ell_t} v_{1} + a_{2,{k_1}+ \ell_t} v_{2} + \cdots + a_{h,{k_1}+ \ell_t} v_{h} = n
\]
has a solution in integers $v_{1},\ldots, v_{h}$ satisfying 
\[
0 \leq v_{i} \leq \frac{(h-1)S_{{k_1}+ \ell_{t+1}}}{a_{i,{k_1}+ \ell_t}}
\]
for $i=1,2,\ldots, h.$  Since $a_{i,k_1+ \ell_t}v_{i} \in \{ 0 \} \cup V(t) \subseteq \mca$, it follows that  $n \in h\mca$ and so $\mca$ is a basis of order $h$.

Next, we prove that the basis \mca\ is thin.  
For all $t \geq 1$ there is a simple upper bound for the cardinality of the finite set $V(t)$:
\[
|V(t)| 
\leq \sum_{i=1}^h \frac{(h-1)S_{k_1+ \ell_{t+1}}}{a_{i,k_1+\ell_t}} 
< \frac{h(h-1)S_{k_1+ \ell_{t+1}}}{a_{h,k_1+\ell_t}}.
\]
Let $A(x)$ be the counting function of the set \mca.  
For $x > (h-1) S_{k_1}$, there is a unique positive integer $t^*$ such that 
\[
(h-1) S_{k_1+\ell_{t^*}} < x \leq (h-1) S_{k_1+\ell_{t^*+1}}.
\]
Using properties of the sequence $\{\ell_t\}_{t=1}^{\infty}$, we obtain 
\begin{align*}
A(x) & \leq A\left( (h-1) S_{k_1+\ell_{t^*+1} } \right)   \\
& \leq (h-1)S_{k_1} 
+  \sum_{t=1}^{t^*}  \frac{ h(h-1)S_{k_1+\ell_{t+1}}}{a_{h,{k_1}+\ell_t}} 
+  \sum_{t =t^*+1}^{\infty} \frac{h(h-1)S_{k_1+ \ell_{t^* +1}}}{a_{h,{k_1}+\ell_t }} \\
& \leq  (h-1)S_{k_1} 
+  \sum_{t=1}^{t^*} \frac{h(h-1)S_{k_1}^{t+1}}{a_{h,{k_1}+\ell_t}}  +  \sum_{t = t^*+1}^{\infty} \frac{h(h-1)S_{k_1}^{t^*+1}}{a_{h,{k_1}+\ell_t}}   \\
& <  h(h-1)S_{k_1} 
\left( 1 +  \sum_{t=1}^{t^*} \frac{S_{k_1}^{t}}{a_{h,{k_1}+\ell_t}}  
+  \sum_{t = t^*+1}^{\infty} \frac{S_{k_1}^{t^*}}{a_{h,{k_1}+\ell_t}} \right)  \\
& \leq  h(h-1)S_{k_1} 
\left( 1 +  \sum_{t=1}^{t^*} \frac{S_{k_1+ \ell_t +1}}{a_{h,{k_1}+\ell_t}}  
+  \sum_{t = t^*+1}^{\infty} \frac{S_{k_1+ \ell_{t^*}+1}}{a_{h,{k_1}+\ell_t}} \right) \\
& \leq  2h(h-1)S_{k_1} 
\left(\sum_{t=0}^{t^*} \frac{S_{k_1+ \ell_t }}{a_{h,{k_1}+\ell_t}}  
+  \sum_{t = t^*+1}^{\infty} \frac{S_{k_1+ \ell_{t^*}}}{a_{h,{k_1}+\ell_t}} \right).  
\end{align*}

We evaluate the sums separately.  Let $t \geq 1$.  
Since $P \geq r_h - r_1$ and $k_1+\ell_t \geq k_1 \geq k_0$, we have 
\[
\frac{s_{h,{k_1}+\ell_t}}{s_{1,{k_1}+\ell_t}} = \frac{({k_1}+\ell_t)P+r_h}{({k_1}+\ell_t)P+r_1} < 1 + \frac{1}{{k_1}+\ell_t} < 2^{1/h}
\]
and so 
\[
s_{h,{k_1}+\ell_t} < 2^{1/h} s_{1,{k_1}+\ell_t}.
\]
For all positive integers $k$, we define 
\[
\sigma_{k}  = S_k^{1/h}.
\]
Then $\sigma_k > 2$ since $S_k > 2^h$.  

The inequality $(h-1)S_{k_1+\ell_{t^*}} < x$ implies that 
\[
\sigma_{k_1+\ell_{t^*}}  \leq (h-1)^{1/h} \ \sigma_{k_1+\ell_{t^*}} 
< x^{1/h}.
\]
We have 
\[
 S_{k_1+\ell_t} = \sigma_{k_1+\ell_t}^h  \leq S_{k_1}^t = \sigma_{k_1}^{ht} < S_{k_1 + \ell_t + 1} < 2S_{k_1 + \ell_t} = 2\sigma_{k_1 + \ell_t} ^h
\]
and so
\[
\sigma_{k_1+\ell_t} \leq \sigma_{k_1}^t 
< 2^{1/h}\ \sigma_{k_1+\ell_t} .
\]
Similarly, the inequality 
\[
s_{1,{k_1}+\ell_t}^h < S_{{k_1}+\ell_t} < s_{h,{k_1}+\ell_t}^h
\]
implies that
\[
s_{1,{k_1}+\ell_t} < \sigma_{k_1+\ell_t}  < s_{h,{k_1}+\ell_t} 
< 2^{1/h} \ s_{1,{k_1}+\ell_t} < 2^{1/h} \  \sigma_{k_1+\ell_t}.
\]
This gives an upper bound for the first sum: 
\begin{align*}
\sum_{t=0}^{t^*} \frac{S_{k_1+ \ell_t }}{a_{h,{k_1}+\ell_t}}  
& = \sum_{t=0}^{t^*} s_{h,{k_1}+\ell_t}  
<  2^{1/h} \sum_{t=0}^{t^*} \sigma_{k_1+\ell_t} \\
& \leq  2^{1/h}\sum_{t=0}^{t^*} \sigma_{k_1}^t 
= 2^{1/h} \left( \frac{  \sigma_{k_1}^{t^*+1}-1}{\sigma_{k_1} -1} \right) \\
& <  2\cdot 2^{1/h} \ \sigma_{k_1}^{t^*} < 2\cdot 4^{1/h} \ \sigma_{k_1+ \ell_{t^*}} \\
& \leq 4 x^{1/h}.
\end{align*}
We used the fact that $S_{k_1} >  2^h$ and so $\sigma_{k_1} > 2$.

To obtain an upper bound for the second sum, we observe that 
\[
\sigma_{k_1}^t < \sigma_{k_1+\ell_t + 1} 
< 2^{1/h} \ \sigma_{k_1+\ell_t} 
\]
and so
\[
\frac{1}{\sigma_{k_1+\ell_t}^{h-1} } 
< \frac{ 2^{(h-1)/h}}{ \sigma_{k_1}^{t(h-1)}}.
\]
Then
\begin{align*}
\sum_{t = t^*+1}^{\infty} \frac{S_{k_1+ \ell_{t^*}}}{a_{h,{k_1}+\ell_t}} 
& = S_{k_1+ \ell_{t^*}}  \sum_{t = t^*+1}^{\infty} \frac{s_{h,k_1+\ell_t}}{S_{k_1+ \ell_t}} \\
& = 2^{1/h} \  S_{k_1+ \ell_{t^*}}
\sum_{t = t^*+1}^{\infty} \frac{\sigma_{k_1+\ell_t}}{\sigma_{k_1+ \ell_t}^h} \\
& <  2^{1/h} \ S_{k_1+ \ell_{t^*}}
\sum_{t = t^*+1}^{\infty} \frac{1}{\sigma_{k_1+ \ell_t}^{h-1}} \\
& <  2^{1/h} \ S_{k_1+ \ell_{t^*}}
\sum_{t = t^*+1}^{\infty} \frac{ 2^{(h-1)/h} }{\sigma_{k_1}^{t(h-1)}} \\
& = 2S_{k_1 + \ell_{t^*}}
\sum_{t = t^*+1}^{\infty} \frac{1}{\sigma_{k_1}^{t(h-1)}} \\
& = \frac{2S_{k_1 + \ell_{t^*}}}{\sigma_{k_1}^{(t^*+1)(h-1)}}
\left(  \frac{\sigma_{k_1}^{h-1}}{  \sigma_{k_1}^{h-1} - 1}  \right)\\
& < \frac{4\sigma_{k_1}^{ht^*}}{\sigma_{k_1}^{(t^*+1)(h-1)}}  \\
& < 4\sigma_{k_1}^{t^*}  \\
& < 8x^{1/h}.
\end{align*}
We conclude that 
\[
A(x) < 24h(h-1)S_{k_1} x^{1/h}
\]
and so \mca\ is a thin basis of order $h$.  
This completes the proof.
\end{proof}

\def\cprime{$'$} \def\cprime{$'$} \def\cprime{$'$}
\providecommand{\bysame}{\leavevmode\hbox to3em{\hrulefill}\thinspace}
\providecommand{\MR}{\relax\ifhmode\unskip\space\fi MR }
\providecommand{\MRhref}[2]{%
  \href{http://www.ams.org/mathscinet-getitem?mr=#1}{#2}
}
\providecommand{\href}[2]{#2}

\end{document}